\RequirePackage{fix-cm}
\documentclass[smallextended]{svjour3}       
\smartqed  

\usepackage{graphicx}
\usepackage[utf8]{inputenc}
\usepackage[english]{babel}
\usepackage{amsmath}
\usepackage{mathtools}
\usepackage{graphicx}
\usepackage{grffile} 
\usepackage{float}
\usepackage{blindtext}
\usepackage{cite}
\usepackage{caption}
\usepackage{subcaption}
\usepackage{color}
\usepackage[normalem]{ulem} 

\newcommand{\squeezeup}{\vspace{-2.5mm}}
\begin{document}

\title{New Numerical Interface Scheme for the Kurganov-Tadmor second-order Method}




\author{Pablo Montes       \and
        Oscar Reula 
}


\institute{M. Pablo $\cdot$ O. Reula \at
              Instituto de Física Enrique Gaviola, CONICET, FaMAF - Universidad Nacional de Córdoba, Ciudad Universitaria, 5000 Córdoba, Argentina \\
              \email{pabloemontes@mi.unc.edu.ar}\\
              \email{oreula@unc.edu.ar}
}

\date{Received: date / Accepted: date}

\maketitle

\begin{abstract}

In this paper, we develop a numerical scheme to handle interfaces across computational domains in multi-block schemes for the approximation of systems of conservation laws. We are interested in transmitting shock discontinuities without lowering the overall precision of the method. We want to accomplish this without using information from interior points of adjacent grids, that is, sharing only information from boundary points of those grids. To achieve this, we choose to work with the second-order Kurganov-Tadmor (KT) method at interior points, relaxing it to first order at interfaces. This allows us to keep second-order overall accuracy (in the relevant norm) and at the same time preserve the TVD property of the original scheme. 
After developing the method we performed several standard one and two-dimensional tests. Among them, we used the one-dimensional advection and Burgers equations to verify the second-order convergence of the method. We also tested the two-dimensional Euler equations with an implosion and a Gresho vortex\cite{liska2003}. In particular, in the two-dimensional implosion test we can see that regardless of the orientation of shocks with respect to the interface, they travel across them without appreciable deformation both in amplitude and front direction.

\end{abstract}


\section{Introduction}

Many physical systems of interest that can be described by partial differential equations (PDEs) cannot be numerically approximated by using a single numerical grid. This may be due to several reasons, for example: a) the topology of the problem does not allow it, as is the case of fields defined on a sphere, b) the size of the problem calls for block-parallelization, 
c) the need for greater resolution in certain parts of the domain. This leads to the creation of methods that allow for multiple grid patches to cover the computational domain. One such type of method is that of multi-block domain decomposition, in which grids are non-overlapping and only connect at boundary regions called interfaces. On those interfaces, the PDE discretization must be specially tailored to allow the correct and accurate transmission of numerical solutions between grids. 

The treatment of these interfaces is generally well understood for smooth solutions. In the finite-difference case satisfying the summation by parts (SBP) property, simultaneous approximation terms (SAT) schemes can be used in a very efficient way and with very high accuracy \cite{CARPENTER1999341}. Alternatively, pseudo-spectral methods are also an optimal choice and similar penalty methods can be applied \cite{HESTHAVEN200023}. However, in systems where shocks or other types of discontinuities are expected, the application  of the above-mentioned schemes is no longer possible. This is mostly due to the formation of spurious oscillations on the numerical approximation of the solution near discontinuities. 

A very large and important subset of systems that present shocks and which are ubiquitous in physics are the so-called systems of conservation laws, which describe the evolution of the density fields of conserved quantities, such as energy, momentum, or mass; that is, fields whose variation satisfies a flux law. Their general form is

\begin{equation}
	\frac{\partial}{\partial t}u + \sum_{j} \frac{\partial}{\partial x_{j}}F^{j}(u) = 0,\quad u|_{\partial \Omega, t=0} = u_{0},
	\label{eq:sysconservationlaws}
\end{equation}

\noindent where $u$ is a vector of conserved quantities, $F^{j}$ are flux functions for each coordinate, $\Omega$ is the computational domain, and $\partial \Omega$ its boundary.
Although these systems present generic shock discontinuities, their structure allows for very good and stable approximations of solutions via so called numerical conservative methods. These are methods that use the weak formulation of Eq. \eqref{eq:sysconservationlaws} to numerically enforce the conservation laws, allowing high orders of convergence in the smooth parts of the solution while keeping good accuracy near shock discontinuities. They also avoid spurious oscillations near shocks by using flux limiters that are built to preserve certain properties of the exact solution. Some of them, for example, preserve monotonicity, others, guarantee that the total variation of the solution does not increase \cite{LeVeque-1992}. To achieve this, conservative methods determine the shock region and locally increase the numerical dissipation.
In the cases of interfaces among multi-blocks, the complication resides in the fact that the information needed to properly determine a shock region and evolve it accordingly is usually located in an adjacent grid. Nevertheless, if we are willing to tolerate a first-order scheme in the interface, the complication is surmountable. Since this region's size is similar to that of a grid cell, the overall global order of convergence can be one order higher.
There are previous works that present ways to treat these interfaces. Their common feature is the use of information from the interior points of adjacent grids; this requires overlapping between grids, thus needing some interpolation and more communication among different processes. Some examples of this can be found in \cite{MCCORQUODALE2015181, Chesshire1994, SHERER2005459, PARTENANDER1994551}.  Here we are interested in developing a simple interface method for conservative schemes which only uses information from the boundary grid points. This is relevant for implementations in high-performance computing codes, for which minimal communication between grids, and by extension computational nodes, is desired. This will also allow the direct treatment of situations where the traversal grid directions are not aligned across interfaces, as could be the case of several grid patches covering a sphere. In this case, it is not straightforward to do a flux reconstruction using interior points from both sides of an interface. 

In this work, we develop a method which overcomes the aforementioned challenges and does comply with the desired properties for an extension of the second-order Kurganov-Tadmor (KT) conservative scheme \cite{KT99}. Our method allows the transmission of shock discontinuities through interfaces while maintaining the global order of convergence and with minimal communication among grids, since only the values at the interface boundary from each grid are shared. Furthermore, it is also built as to preserve the total variation diminishing (TVD) property of the original semi-discrete KT scheme. Contrary to upwind methods, KT does not require a characteristic decomposition of the fields, making it easier to implement.

This paper is organized as follows: In section \ref{sec:KT-method}, we give a brief overview of the KT method and its total variation diminishing (TVD) property. In section \ref{sec:interface-method}, we introduce our multi-block version of the KT method and show that the TVD property is still maintained. It works by falling back to a specific first-order formulation near the boundary of each grid, and later calculating a cell average between the values obtained in adjacent grids. Finally, in section \ref{sec:results}, we apply our scheme to the case of the advection and Burgers equations in one dimension to test convergence and shock propagation. We also corroborate that the second-order convergence is maintained and shock discontinuities and waves are completely transmitted through interfaces without loss of neither phase nor group speed. Finally, we show results for the implementation of this interface for the two-dimensional Euler equations, solving a Gresho vortex and an implosion problem from \cite{liska2003}.

\section{Kurganov-Tadmor method}
\label{sec:KT-method}

In this section we make a quick overview of the semi-discrete formulation of the KT method. More information can be found in \cite{KT99}.

We start with a simple one-dimensional problem of the form

\begin{equation}
    \frac{\partial u}{\partial t}(x,t) + \frac{\partial f(u(x,t))}{\partial x}=0 \quad (x,t) \in[0,\infty)\times[0,\infty),
    \label{eq:scalarconservationlaw}
\end{equation}

\noindent where $u$ is a vector of conserved quantities and $f$ are the corresponding fluxes.

We start by discretizing the spatial dimension in equal intervals of width $\Delta x$, and denoting $x_j \doteq j\Delta x$. We then define the cell average of $u$ at grid point $j$  as

\begin{equation}
    \overline{u}_j(t) = \frac{1}{\Delta x}\int_{x_{(j-\frac{1}{2})}}^{x_{(j+\frac{1}{2})}}u(\xi,t)\;d\xi.
    \label{eq:cellaverage}
\end{equation}

If we now combine \eqref{eq:scalarconservationlaw} and \eqref{eq:cellaverage} we obtain an equivalent weak formulation of the conservation law,

\begin{equation}
    \frac{\partial \overline{u}_{j}}{\partial t}(x,t) =  f\left ( u(x_{j-1/2,},t)\right) - f\left ( u(x_{j+1/2,},t)\right).
    \label{eq:weakscalarconservationlaw}
\end{equation}

Conservative methods evolve these cell averages of the solution using Eq. \eqref{eq:weakscalarconservationlaw} and thus need to approximate the solution at the cell boundaries $x_{j-1/2}$ and $x_{j+1/2}$. The main difference between different schemes lies in the way that these boundaries values are approximated. From now onwards we denote $v_j(t)$ as the numerical approximation of  $\overline{u}(x_j, t)$. 


The semidiscrete KT method takes the conservative form

\begin{equation}
    \frac{\mathrm{d}}{\mathrm{d}t}v_j(t) = -\frac{H_{j+1/2}(t)-H_{j-1/2}(t)}{\Delta x},
\end{equation}

\begin{eqnarray}
    H_{j+1/2}(t):=\frac{f(v^+_{j+1/2}(t))+f(v^-_{j+1/2}(t))}{2}-
\frac{a_{j+1/2}(t)}{2}\left[v^+_{j+1/2}(t)-v^-_{j+1/2}(t) \right],
\end{eqnarray}

\begin{equation}
    v^{+}_{j+1/2}(t) = v_{j+1}(t) - \frac{\Delta x}{2} (v_x)_{j+1}(t) \quad v^{-}_{j+1/2}(t) = v_{j}(t) + \frac{\Delta x}{2} (v_x)_{j}(t),
    \label{eq:upumdefinition}
\end{equation}

\noindent where $(v_x)_j$ is some approximation of the spatial derivative of $v$, $a_{j+1/2}(t)$ the maximum propagation speed of the solution between $v^-_{j+1/2}(t)$ and $v^+_{j+1/2}(t)$,

\begin{equation}
    a_{j+1/2}(t) = \max \left (  \rho \left (\left |\frac{\partial f}{\partial v} \left (v^+_{j+1/2}(t) \right) \right |\right ), \rho \left (\left |\frac{\partial f}{\partial v} \left (v^-_{j+1/2}(t)\right )  \right | \right )\right ),
    \label{eq:adefinition}
\end{equation}

\noindent and $\rho (|\frac{\partial f}{\partial v}(v)|)$ represents the spectral radius of the matrix $|\frac{\partial f}{\partial v}(v)|$.

While this scheme may look complicated at first glance, a good starting point to understand it is to set $(v_x)_j = 0$. This assumption gives us the first order KT method, which can be expressed in the dissipative form

\begin{equation}
    \begin{split}
	\frac{d}{dt} (v_j) = & -\frac{1}{2\Delta x}\left( f(v_{j+1})-f(v_{j-1})\right )+\\&{\frac{1}{2\Delta x} \left(a_{j+1/2}(v_{j+1} - v_{j})-a_{j-1/2}(v_{j} - v_{j-1})\right)}.
    \end{split}
    \label{eq:KT-dissipative-form}
\end{equation}

While the first term in the right-hand side of Eq. \ref{eq:KT-dissipative-form} corresponds to the approximation of the numerical flux, the second one corresponds to a numerical dissipation, being similar to a numerical approximation of the second derivative of $u$. Because $(v_x)_j = 0$, we can interpret our numerical approximation as being piecewise constant in the intervals $x_{j-1/2}, x_{j+1/2}$. 

In contrast, the second-order KT method does not set $v_x$ to zero and is piecewise \textit{linear} in those same intervals. The approximation of $v_x$ is chosen as to preserve the \textit{total variation diminishing} property (TVD). Recall that the \textit{total-variation} (TV) in the case of a numerical solution is defined as 

\begin{equation}
    TV(v(t)) = \sum_{j}|v_{j+1}(t)-v_{j}(t)|,
\end{equation} 

\noindent while in the case of an exact solution, the TV is defined as 
\begin{equation}
 TV(u(t)) = \int_\Omega\left |\frac{\partial u}{\partial{x}}(x,t) \right |dx.   
\end{equation}
The diminishing of the TV quantity is called the TVD property. In the case that the quantity is differentiable, it implies that $\frac{d}{dt} TV(v) \leq 0$. This property holds for exact solutions of scalar one-dimensional conservation laws and it ensures that no new extrema are created in the evolution of the numerical solution, thus preventing spurious oscillations near shock discontinuities\cite{LeVeque-1992}. Note that so far we are referring to the semi-discrete formulation of the KT scheme; to preserve the TVD property for the fully discretized system we also \sout{need to} require the use of special temporal integrators \cite{tvd-runge-kutta}.

The following choice for $v_x$ satisfies the TVD property \cite{TVDproperty}.

\begin{equation}
    (v_x)_j = \textrm{minmod} \left(\theta \frac{v_{j}-v_{j-1}}{\Delta x}, \frac{v_{j+1}-v_{j-1}}{2\Delta x}, \theta \frac{v_{j+1}-v_{j}}{\Delta x}\right ),
    \label{eq:minmodderivative}
\end{equation}

\noindent with the function $\textrm{minmod}$ defined as

\begin{equation}
\textrm{minmod} (x_1, x_2, x_3,...) = \left\{ \begin{matrix}
 \max( x_{j} ) & \textrm{if } x_{j}<0 \;\forall j,\\
 \min( x_{j} ) & \textrm{if } x_{j}>0 \;\forall j,\\ 
 0 & \textrm{otherwise,}
\end{matrix} \right.
    \label{eq:minmoddef}
\end{equation}
and $\theta$ a parameter between $1$ and $2$ \cite{KT99}.

For systems of equations, the TVD property is not valid. For instance, consider the Euler equations and initial data with zero momentum density in a compact region. As the system starts to evolve, the fluid starts to move in different directions, and so the total variation of the momentum density increases. 
Nevertheless, enforcing this property for the scalar case singles out a particular one-parameter family of choices of derivatives. That choice is then enforced even for systems of equations and even in several dimensions, where, in practice, shows very nice stability properties and stops the appearance of local oscillations near shocks.  
These methods converge to second-order in the $Lip'$ semi-norm  \cite{KT99}, also known as the Wasserstein distance, which is especially suited for analyzing the convergence of conservative methods. The convenience of this norm for weak solutions is very well explained by U. S. Fjordholm and S. Solem in  \cite{sec-ord-conv}. Therefore, this is the  appropriate norm to test the convergence of our interface method. 


\section{Kurganov-Tadmor interface generalization}
\label{sec:interface-method}
\indent In this section, we explain our interface method. For simplicity, we start with the one-dimensional case. First, assume that we choose to place our interface at a point $x_{I}$ of a domain $\Omega$. We label our subdomains to the left and right of $x_{I}$ with the letters $L$ and $R$ respectively.

\begin{figure}[H]
    \centering
    \includegraphics[width = 0.9\linewidth]{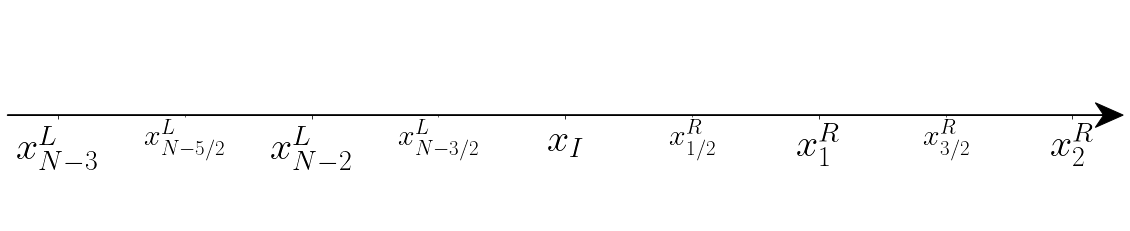}
    \caption{We split our domains in two regions $L$ and $R$ around an interface $x_{I}$, and discretize both regions so that they share a grid point at the interface. We mark our intermediate points for better understanding of the scheme. Even though in the figure $\Delta x^{L} = \Delta x^{R}$, this is not a necessary condition.
}
    \label{fig:interface-line}
\end{figure}

We discretize $L$ using $N$ points and $R$ using $M$ points, with $x_{N-1}^{L} = x_{I} = x_{0}^{R}$, as seen in figure \ref{fig:interface-line}; that is, in both grids we store information of our estimation of $\overline{u}
(x_{I}, t)$.

In other words, we are covering our one-dimensional domain with two grid patches that do not overlap each other and share a single interface point $x_I$. For simplicity, given any variable $w$, we will use the notation $w^{L}_{I-1/2} \doteq w^{L}_{N-3/2}$ and $w^{R}_{I+1/2} \doteq w^{R}_{1/2}$. We only focus on how to approximate fields in the interface, as the rest of the grid is solved using the standard second-order Kurganov-Tadmor method previously explained.

We recall that KT updates the estimated point averages of the solution at each grid point. We want to estimate the cell average of the solution in the interval $(x_{I-1/2}^L,x_{I+1/2}^R)$. The idea of our interface method is to use the points of the subgrids $L$ and $R$ to estimate the left and right averages $v^{L}_{I}(t)$ and $v^{R}_{I}(t)$,

\begin{equation}
v^{L}_{I}(t) \approx \frac{1}{\Delta x^{L}}\int_{x^{L}_{I-1/2}}^{x_{I}} u(\xi, t)\;d\xi, \quad v^{R}_{I}(t) \approx \frac{1}{\Delta x^{R}}\int_{x_{I}}^{x^{R}_{I+1/2}} u(\xi, t)\;d\xi.
\end{equation}

\noindent and latter communicate both grids by calculating the full average $u_{I}$,

\begin{equation}
    v_{I}(t) \doteq \frac{\Delta x^{L} v^{L}_{I}(t) + \Delta x^{R} v^{R}_{I}(t)}{\Delta x^{L} + \Delta x^{R}} \approx \frac{1}{\Delta x^{L} + \Delta x^{R}}\int_{x^{L}_{I-1/2}}^{x^{R}_{I+1/2}} u(\xi, t)\;d\xi.
\end{equation}

\begin{figure}
    \centering
    \includegraphics[width = 0.9\linewidth]{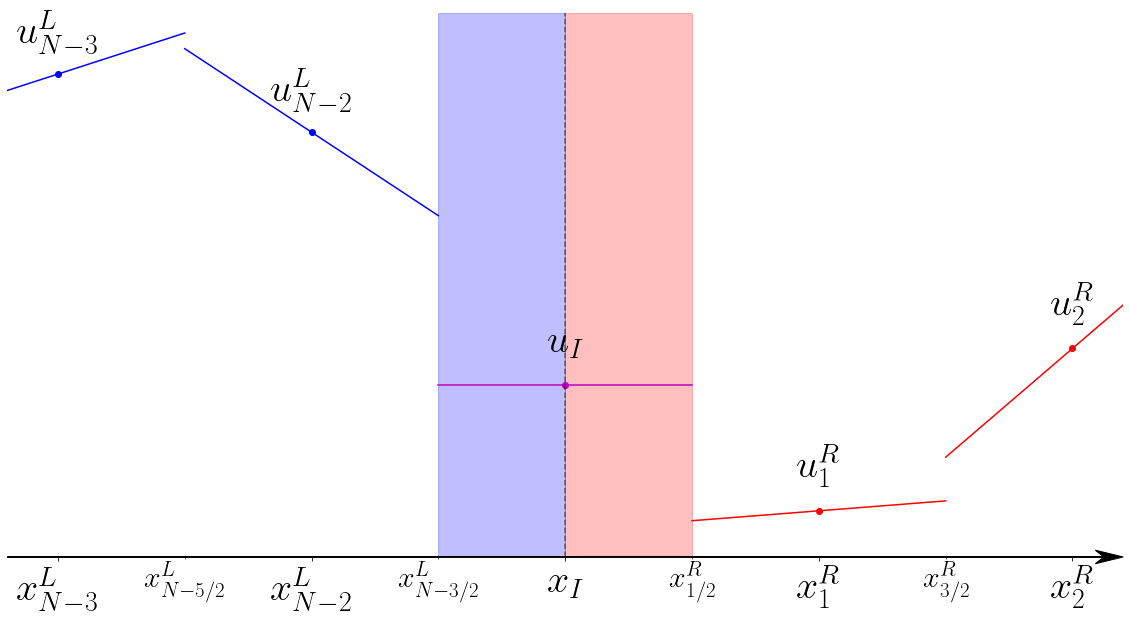}
    \caption{The interface is marked with a dashed line. For each point, we also plot a line of slope $v_x$. The left and right extremes of those lines correspond to the values $v^{+}_{j-1/2}$ and $v^{-}_{j+1/2}$, respectively. We mark in blue and red the averaging areas of $v^{L}_{I}$ and $v^{R}_{I}$.}
    \label{fig:interface-diagram}
\end{figure}

Our method then proceeds as follows. At points near the interface, we lower the approximation to first order. Thus we don't need to compute the derivative at the interface, as seen in Figure \ref{fig:interface-diagram}. This way we avoid using information from the adjacent grid, as normally we would need the values $v_{I-1/2}^L$ and $v_{I+1/2}^R$ to calculate this derivative. An order reduction is in any case expected for touching grids methods, as happens with SAT methods \cite{CARPENTER1999341}. In SAT methods this is due to the fact that SBP operators are lower-order near the interfaces.

At this point, we would still need the values $v_{I-1}^n$ and $v_{I+1}^n$ to calculate the fluxes to both sides of $x_I$. To avoid using the value at the other side of the interface, we instead choose to calculate partial cell averages on a smaller cell, for instance, for the left to the interface we consider the cell between $x_{I-1/2} = x_{N-3/2}$ and $x_I$,  of width $\Delta	x^{L} / 2$. Then, to calculate

\begin{equation}
(v_{I}^L)^{n+1} \approx \int_{x^L_{I-1/2}}^{x^I} u(x,t^{(n+1)})dx,
\end{equation}

we use the reconstructed fluxes at $x^{L}_{I-1/2}$ and $x_{I}$. The algorithm for the left grid becomes

\begin{equation}
\dfrac{d}{dt}v^{L}_{I} = -\dfrac{G_{I}^{L+} - G_{I}^{L-} }{\Delta x^{L}},
\label{eq:interfacestep1}
\end{equation}
with 
\begin{equation}
\begin{matrix}
G_{I}^{L+} &=&   f(v^{L+}_{I-1/2}) + a^{L}_{I-1/2}v^{L+}_{I-1/2}\\
&=&f(v^{L+}_{I-1/2}) + a^{L}_{I-1/2}v_{I},\quad\;\;\;\;\;
\end{matrix} 
\end{equation}
and
\begin{equation}
    G_{I}^{L-} = f(v^{L-}_{I-1/2}) + a^{L}_{I-1/2}v^{L-}_{I-1/2},
\end{equation}


\noindent where the values $v^{L\pm}_{I-1/2}$ and $a^{L}_{I-1/2}$ are calculated using the definitions \eqref{eq:upumdefinition} and \eqref{eq:adefinition}. Because $(v_x)^{L}_{I} = 0$ we also have $v^{L+}_{I-1/2} = v_{I}$.

Similarly, the algorithm for the right grid becomes

\begin{equation}
\dfrac{d}{dt}v^{R}_{I} =-\dfrac{G_{I}^{R+} - G_{I}^{R-} }{\Delta x^{R}},
\label{eq:interfacestep1-right}
\end{equation}

\noindent with
\begin{equation}
G_{I}^{R+} = f(v^{R+}_{I+1/2}) - a^{R}_{I+1/2}v^{R+}_{I+1/2}
\end{equation}
and
\begin{equation}
\begin{matrix}
G_{I}^{R-} &=&   f(v^{R-}_{I+1/2}) - a^{R}_{I+1/2}v^{R-}_{I+1/2},\\
&=&f(v^{R-}_{I+1/2}) - a^{R}_{I+1/2}v_{I}.\quad\;\;\;\;\;
\end{matrix} 
\end{equation}
Each of these calculations is done with information present on their respective side of the interface. Thus, they can be performed without any communication from the other grid. After obtaining both temporal derivatives, we integrate in time using a TVD temporal scheme to evolve $v^{L}_{I}$ and $v^{R}_{I}$. Finally, we communicate both grids by calculating a full grid average between $x_{I-1/2}$ and $x_{I+1/2}$ by averaging $v^{L}_I$ and $v^R_I$:

\begin{equation}
    \begin{array}{r c l}
    v_{I} &=&  \dfrac{\Delta x^{L}v^{L}_{I}+\Delta x^{R}v^{R}_{I}}{\Delta x^{L} + \Delta x^{R}} \\ 
    v^{L}_{I} &\leftarrow & v^{\phantom{R}}_{I}\\
    v^{R}_{I} &\leftarrow & v^{\phantom{L} }_{I}.
\end{array}
\end{equation}

We notice that the propagation information of the waves is only passed at this very last instance, namely when we average the new values on each side using the new values from the other.
This way we essentially recover the original first-order Kurganov-Tadmor method.

We can check that the TVD property still holds. From \cite{KT99} and \cite{TVDproperty} we know that this property is guaranteed if

\begin{equation}
    \left |\frac{\Delta x(v_x)_{j+1/2}}{\Delta v_{j\pm 1/2}}
    \right| \leq 2,
    \label{eq:tvd-condition}
\end{equation}

\noindent with $\Delta v_{j+1/2} = v_{j+1} - v_{j}$. This is clearly the case at the interface point, since $(v_x)_I$ = 0.
Although we could, in theory, approximate $(v_x)^{L}_{I} = \frac{v_{I} - v^{L}_{N-2}}{\Delta x^{L}}$ and $(v_x)^{R}_{I} = \frac{v^{R}_{1}-v_{I}}{\Delta x^{R}}$, doing so would not guarantee that \eqref{eq:tvd-condition} holds.

As a remark, we note that $\dfrac{d v^{L}_{I-1}}{dt}$ and $\dfrac{d v^{R}_{I+1}}{dt}$ are also slightly modified by our method. For example, to calculate $\dfrac{d v^{L}_{I-1}}{dt}$, we need to use equation \eqref{eq:upumdefinition}, which in turns uses the value $(v_x)^{\phantom{R}}_{I}$ which is set to zero at the interface. This is only a slight modification to the standard second-order KT method, thus not altering the TVD condition by the same arguments used before.

As is the case with the full KT scheme, the extension of this method to equation systems and several dimensions is done by implementing it to each vector component and each dimension. We simply need to set the gradient in the normal direction to zero for the interface values, evolve the truncated KT scheme for each side, and then calculate a final cell average using the volumes of the cells in each grid.

\section{Numerical tests}
\label{sec:results}

We present convergence results for one-dimensional scalar equations, as well as results for the two-dimensional Euler equation. In all cases we use the third order TVD Runge-Kutta scheme introduced in equation $(3.3)$ in \cite{tvd-runge-kutta}. We set the parameter $\theta=2.0$ in the minmod definition \eqref{eq:minmodderivative} in the scalar cases, and set $\theta = 1.2$ for the two dimensional Euler equation (see discussion of optimal $\theta$ values in \cite{KT99}). The CFL condition is set so that $\Delta t = 0.1\Delta x$ for all cases.

\subsection{Advection}

We consider the one-dimensional advection equation with speed 1,
\begin{equation}
    \frac{\partial}{\partial t} u(x,t) + \frac{\partial}{\partial x}u(x,t) = 0,
    \label{eq:advection-eq}
\end{equation}
\noindent in the domain $x\in[-2,2)$, with periodic boundary conditions. The initial data is given by a $C^{7}$ peak function modeled with an eighth order polynomial,

\begin{equation}
u_{0}(x)=\left\{\begin{matrix}
 (x-1)^4(x+1)^4&  \textrm{if } |x| < 1\\ 
0 & \textrm{if } |x| \geq 1.  
\end{matrix}\right.    
\label{eq:advection-inidat}
\end{equation}
The exact solution to this problem is trivial and corresponds to $u(x,t) = u_{0}(x-t)$. Since the solution is periodic, we also have $u(x,4n) = u_0(x)  \; \forall n \in \mathcal{Z}$. 

\begin{figure}[H]
    \centering
    \begin{subfigure}[t]{0.48\textwidth}
        \centering
        \includegraphics[width=\textwidth]{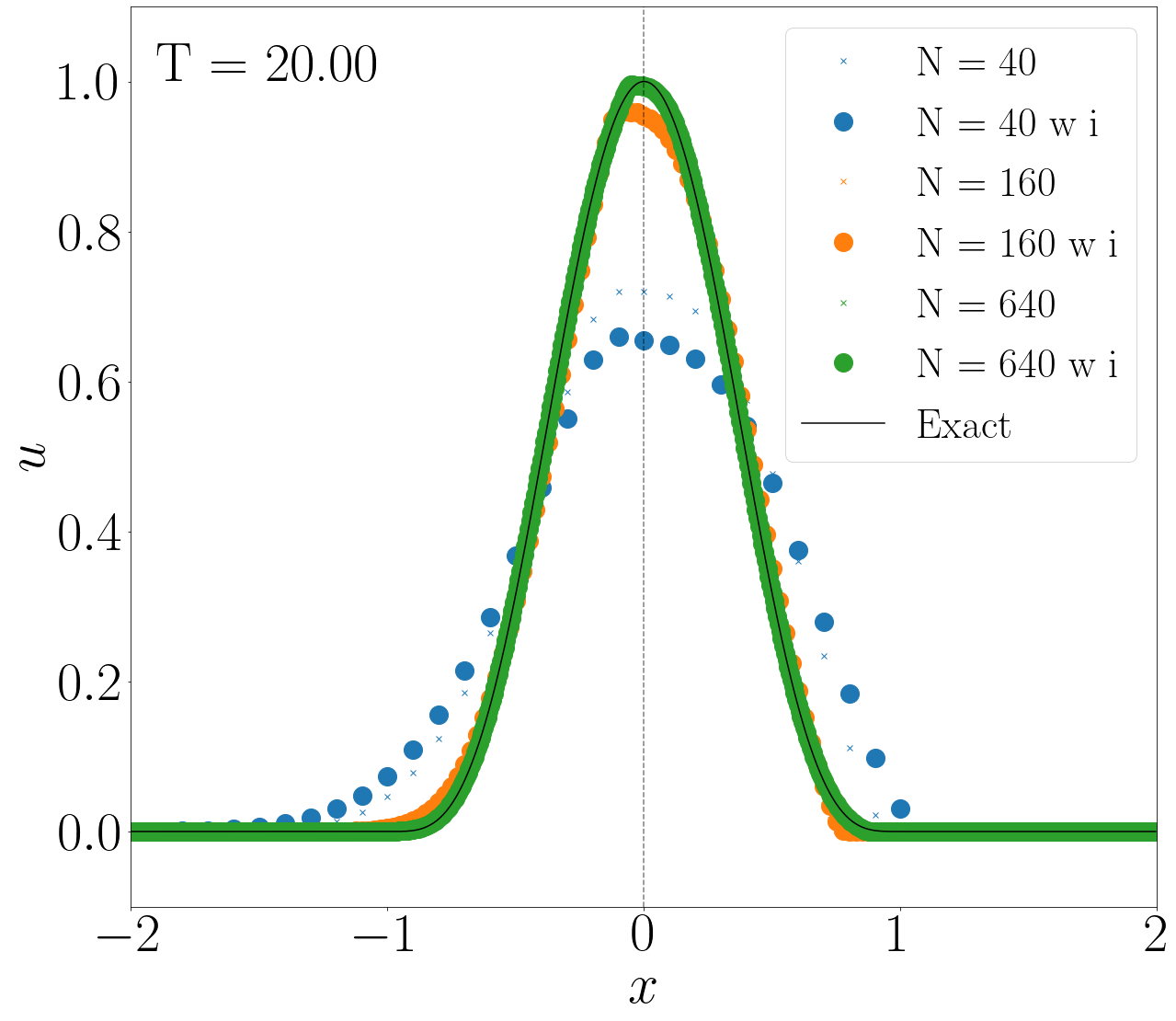}
        \caption{Evolution of \eqref{eq:advection-inidat} for the advection equation with different number of points $N$ and with (w i) and without an interface in $x = 0$.}
        \label{fig:chichon-evolution}
    \end{subfigure}
    \hfill
    \begin{subfigure}[t]{0.48\textwidth}
        \centering
        \includegraphics[width=\textwidth]{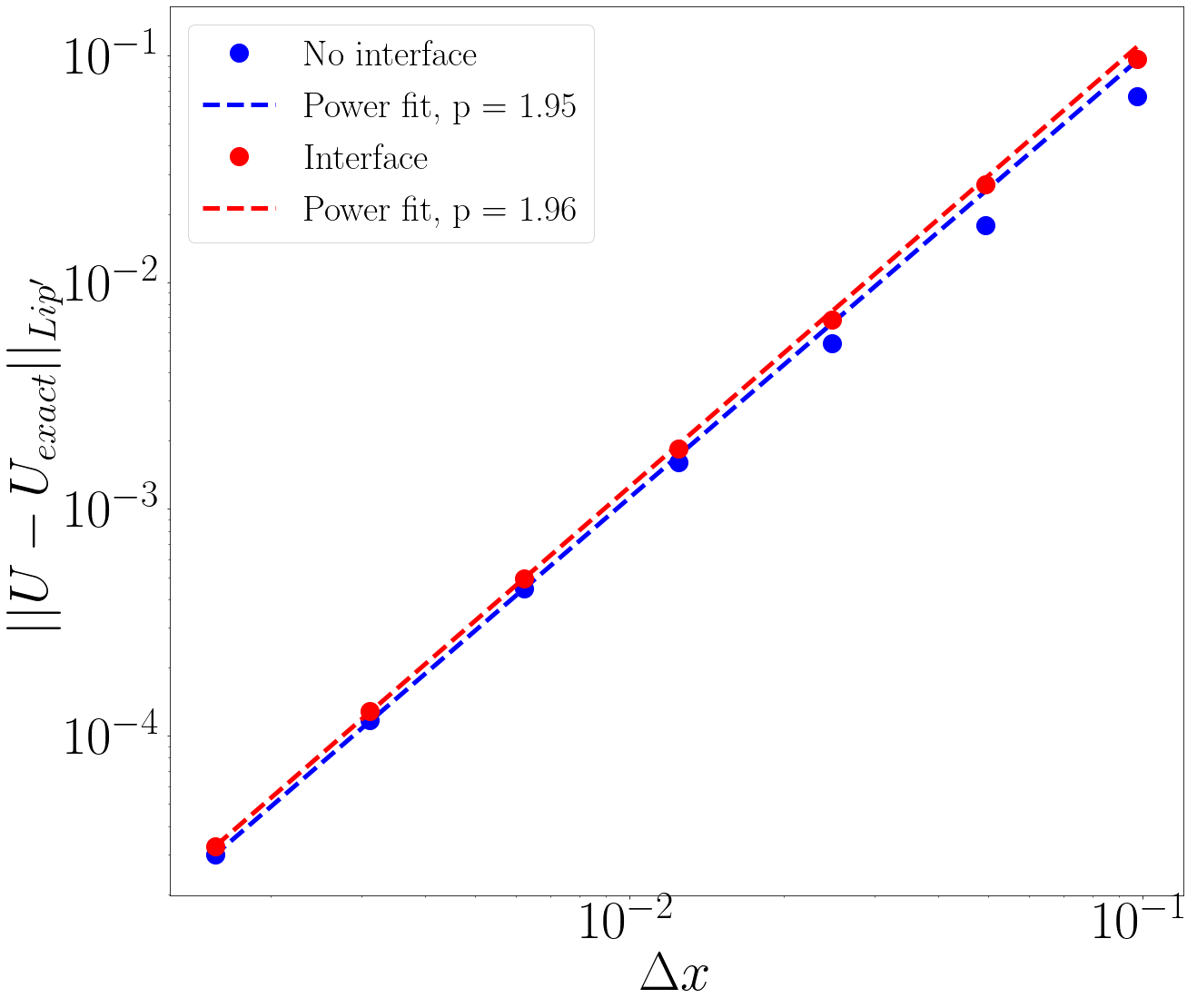}
        \caption{$Lip'$ convergence of \eqref{eq:advection-inidat} for the advection equation.}
        \label{fig:chichon-convergence}
    \end{subfigure}
    \caption{Results for the one-dimensional advection problem at T = 20, corresponding to 4 periods. The CFL is set to $0.1$ in all cases.}
\end{figure}

We evolve \eqref{eq:advection-eq}-\eqref{eq:advection-inidat} until $T = 20$, both with and without a numerical interface placed at $x_{I} = 0.5$. For simplicity in the convergence analysis, we assume $\Delta x^{L} = \Delta x^{R}$ in the interface simulations. The evolution results are shown in Figure \ref{fig:chichon-evolution}, and convergence results are shown in \ref{fig:chichon-convergence} for the seminorm $Lip'$. We fit the convergence data with a power equation

\begin{equation}
    ||u-u_{exact}||_{Lip'} = b\Delta x^p.
    \label{eq:powerfit}
\end{equation}

We notice that even though there is slightly greater numerical dissipation in the presence of interfaces, as seen in the deformation of the initial profile in Figure \ref{fig:chichon-evolution}, the second-order global convergence of the method is maintained. The speed of the wave is also correctly preserved, as evidenced by the approximations still being centered after 5 periods.

\begin{figure}[H]
        \centering
        \includegraphics[width=0.8\textwidth]{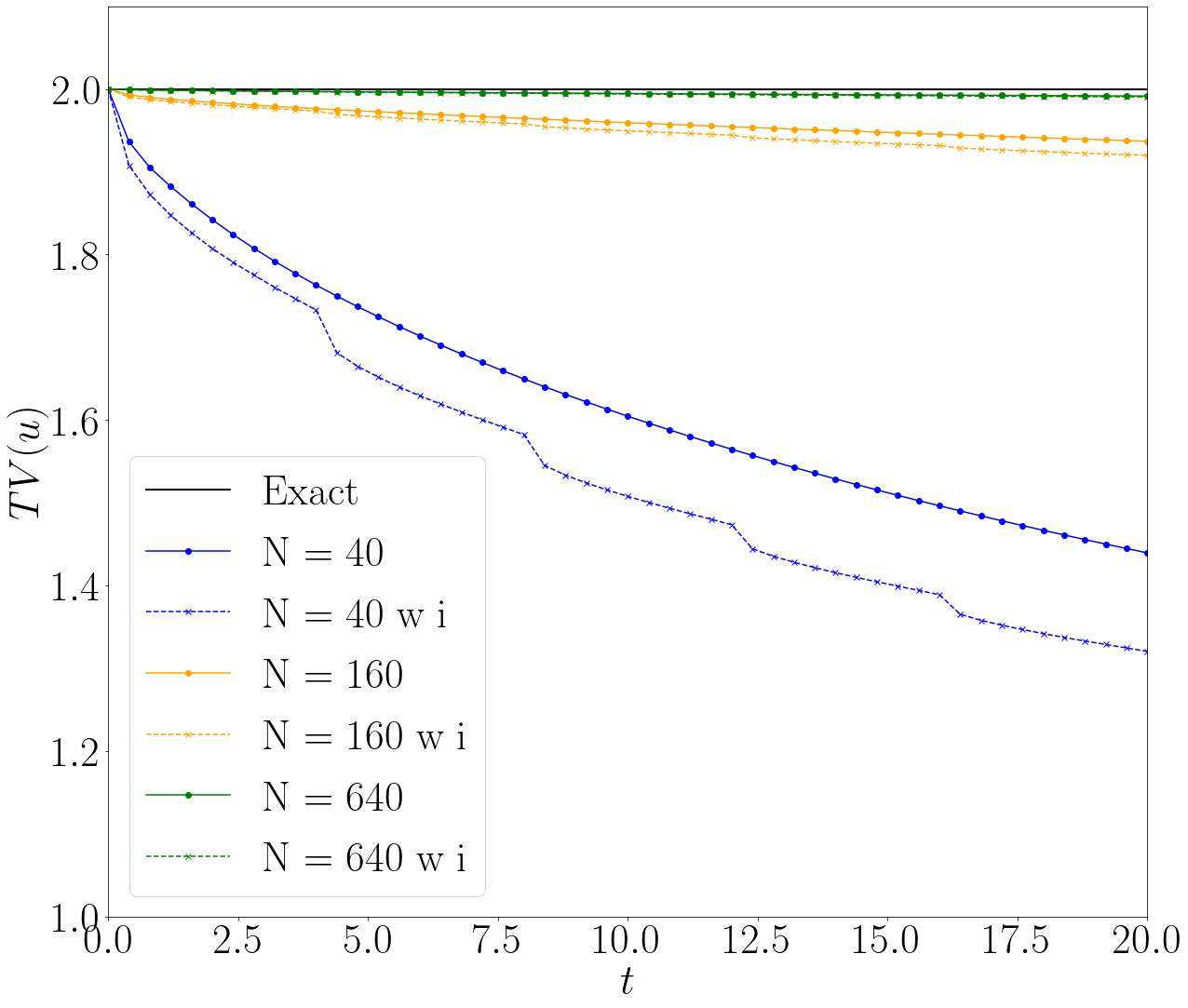}
    \caption{Evolution of the total variation for the advection equation with (w i) and without an interface. The total variation is constant until the shock is developed.}
    \label{fig:tv-advection-chichon}
\end{figure}

Figure \ref{fig:tv-advection-chichon} shows the total variation evolution of the solution of the advection equation until $T = 40$ for several grid resolutions. We notice that the total variation of the numerical solution diminishes, even though the total variation of the exact solution remains constant. We also notice that the total variation diminishes faster in the interface simulations, particularly so whenever the peak of the function passes through the interface. This is due to the fact that the total variation only depends on the extrema of the function, and is therefore in this case only sensitive to the dissipation of the peak of the numerical solution. It is well known that the KT schemes are very dissipative at extrema, thus when a maximum or minimum passes through the interface, where the scheme is only first order, a substantial diminishing of the TV is expected. We shall see in the Burgers equation example that in the case of a shock the situation is much better.

\subsection{Burgers equation}

To study the behavior of shocks passing through an interface we consider now the one-dimensional Burgers equation

\begin{equation}
    u_{t} + \left(\frac{u^{2}}{2}\right)_x = 0,
    \label{eq:Burgers}
\end{equation}
with smooth periodic initial data in the interval $[0, 2\pi)$,

\begin{equation}
    u(x,0) = 0.5 + \sin(x).
    \label{eq:sininidat}
\end{equation}

The solution of \eqref{eq:sininidat} is well known and develops a shock at time $T=1$ with speed $0.5$. Results of the second-order convergence for the KT scheme for this problem in the case without interfaces can be found in \cite{KT99}.

\begin{figure}[H]
        \centering
        \includegraphics[width=0.5\textwidth]{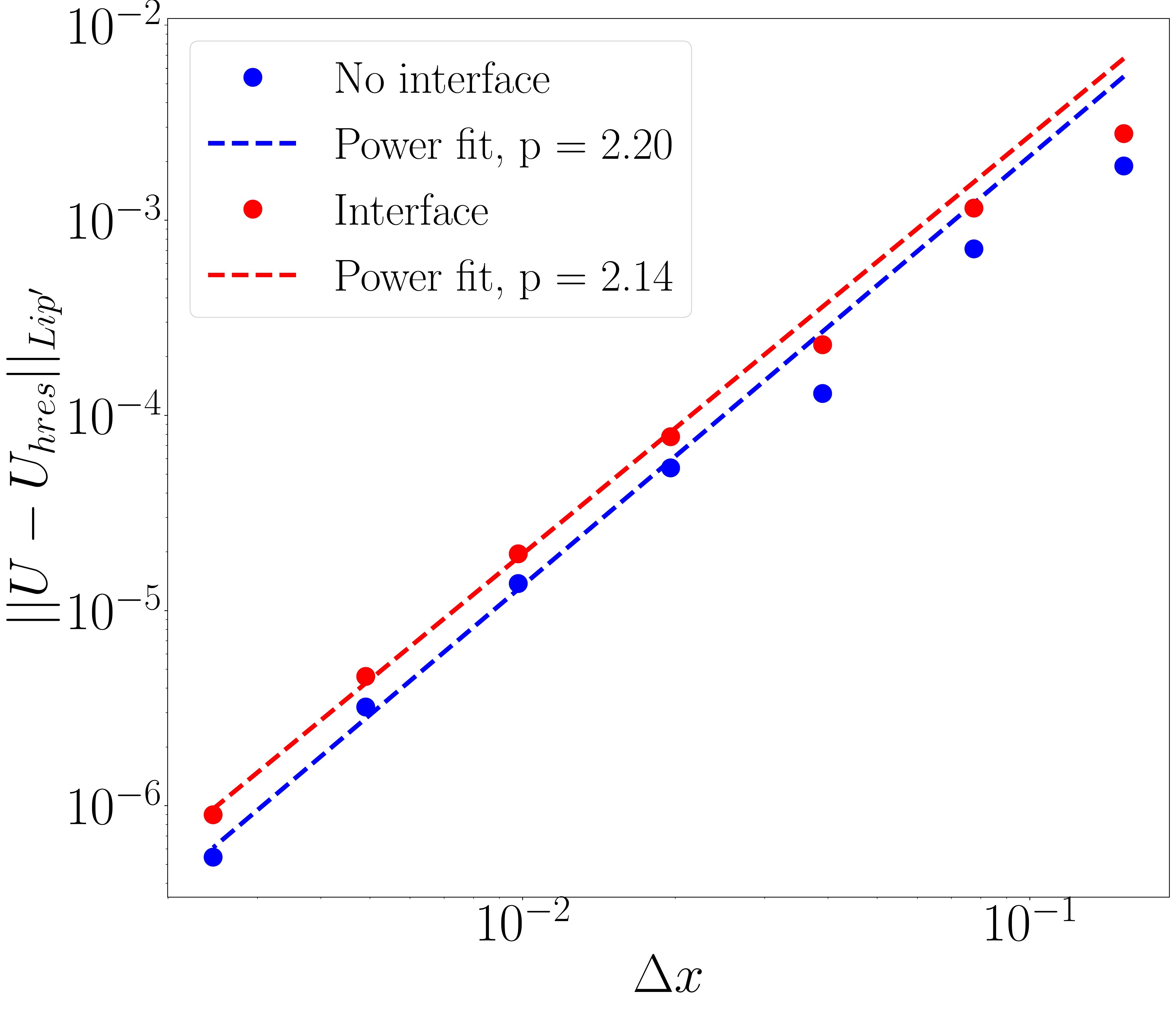}
    \caption{Convergence of \eqref{eq:sininidat} with the Burgers equation at time $T=2.0$ with CFL = $0.1$. $u_{hres}$ corresponds to the approximation with $N = 20480$.}
    \label{fig:sinconvergence1}
\end{figure}

We first evolve the wave until $T = 2.0$ so that the shock has time to develop and travel across the interface, which is placed at $x = 1.25\pi$. The convergence results are shown in Figure \ref{fig:sinconvergence1}. We compare our results with a high resolution approximation using $N=20480$ points without an interface. We notice that the shock discontinuity successfully travels across the interface and the second-order global convergence is maintained .

\begin{figure}[]
        \centering
        \includegraphics[width=0.8\textwidth]{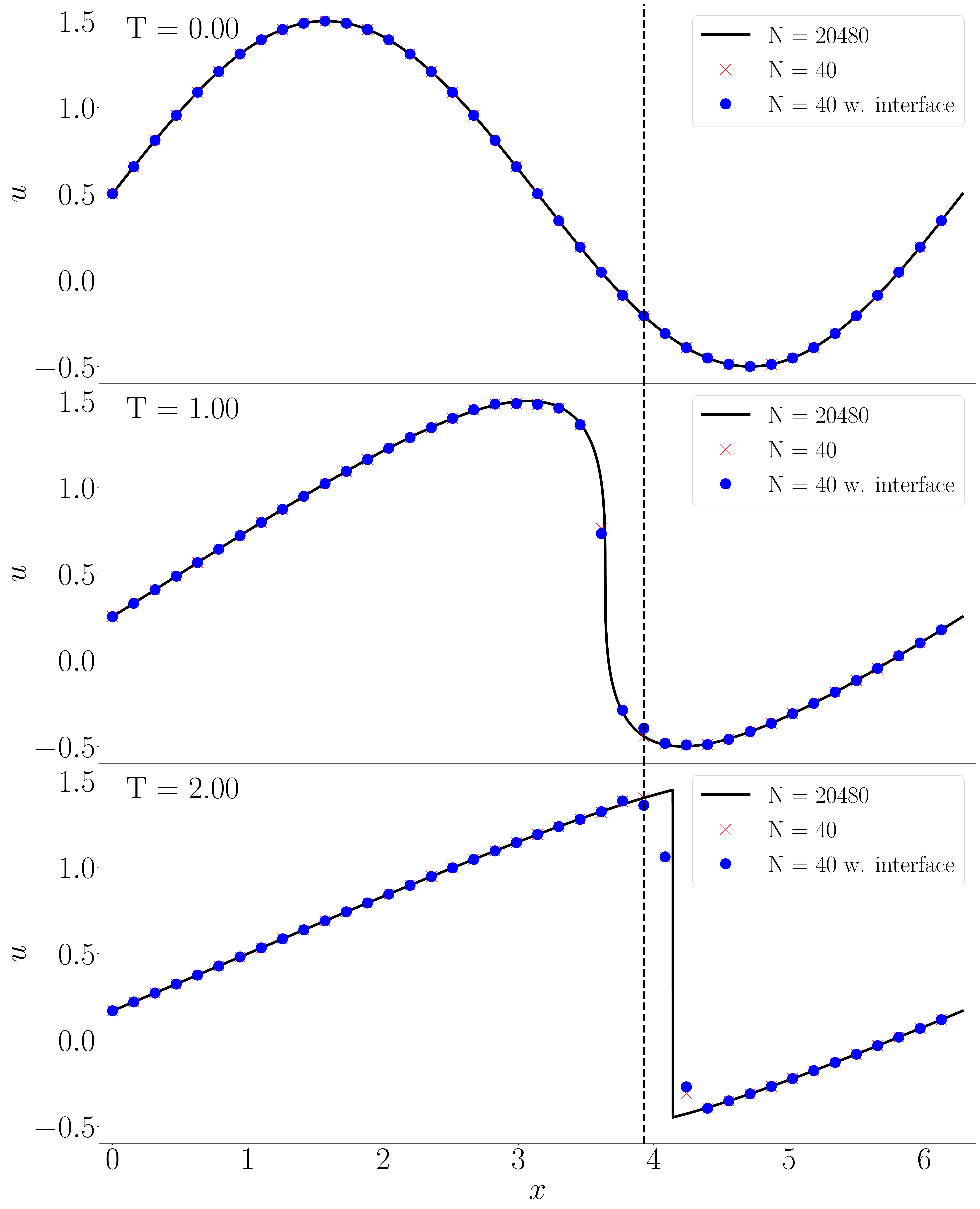}
        \caption{Evolution of \eqref{eq:sininidat} with the Burgers equation at different times for $N = 40$ and $CFL = 0.1$ with and without an interface. The interface location is marked by the dashed black line.}
        \label{fig:sinevolutionshort}

\end{figure}

\begin{figure}[H]
    \centering
    \begin{subfigure}[t]{0.48\textwidth}
        \centering
        \includegraphics[width=\textwidth]{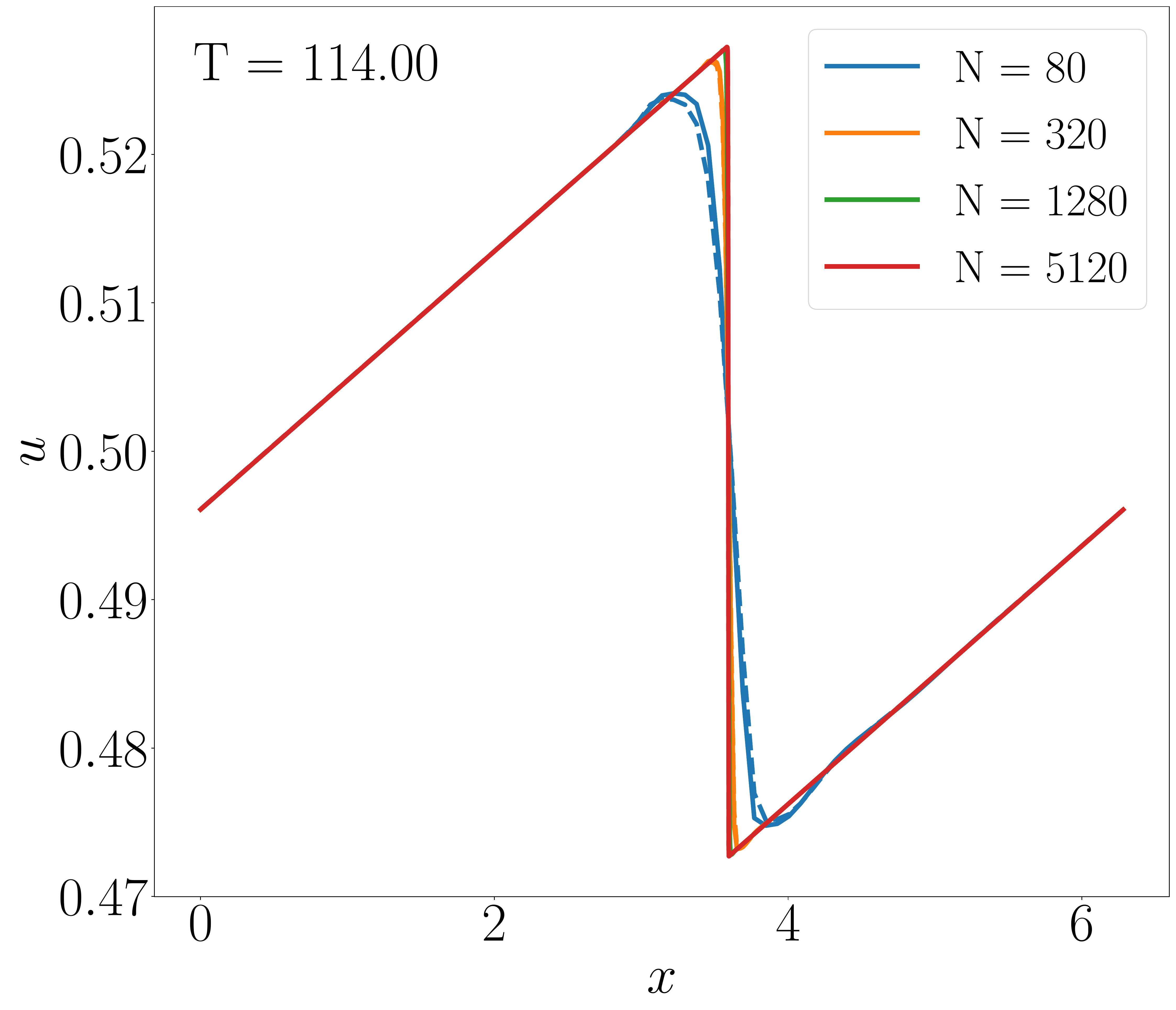}
    \caption{Evolution of \eqref{eq:sininidat} with the Burgers equation with different numbers of points $N$ and with (dashed line) and without (continuous line) interface. This is after 10 passes of the shock across the interface.}
    \label{fig:sinevolutionlong}

    \end{subfigure}
    \hfill
    \begin{subfigure}[t]{0.48\textwidth}
        \centering
        \includegraphics[width=\textwidth]{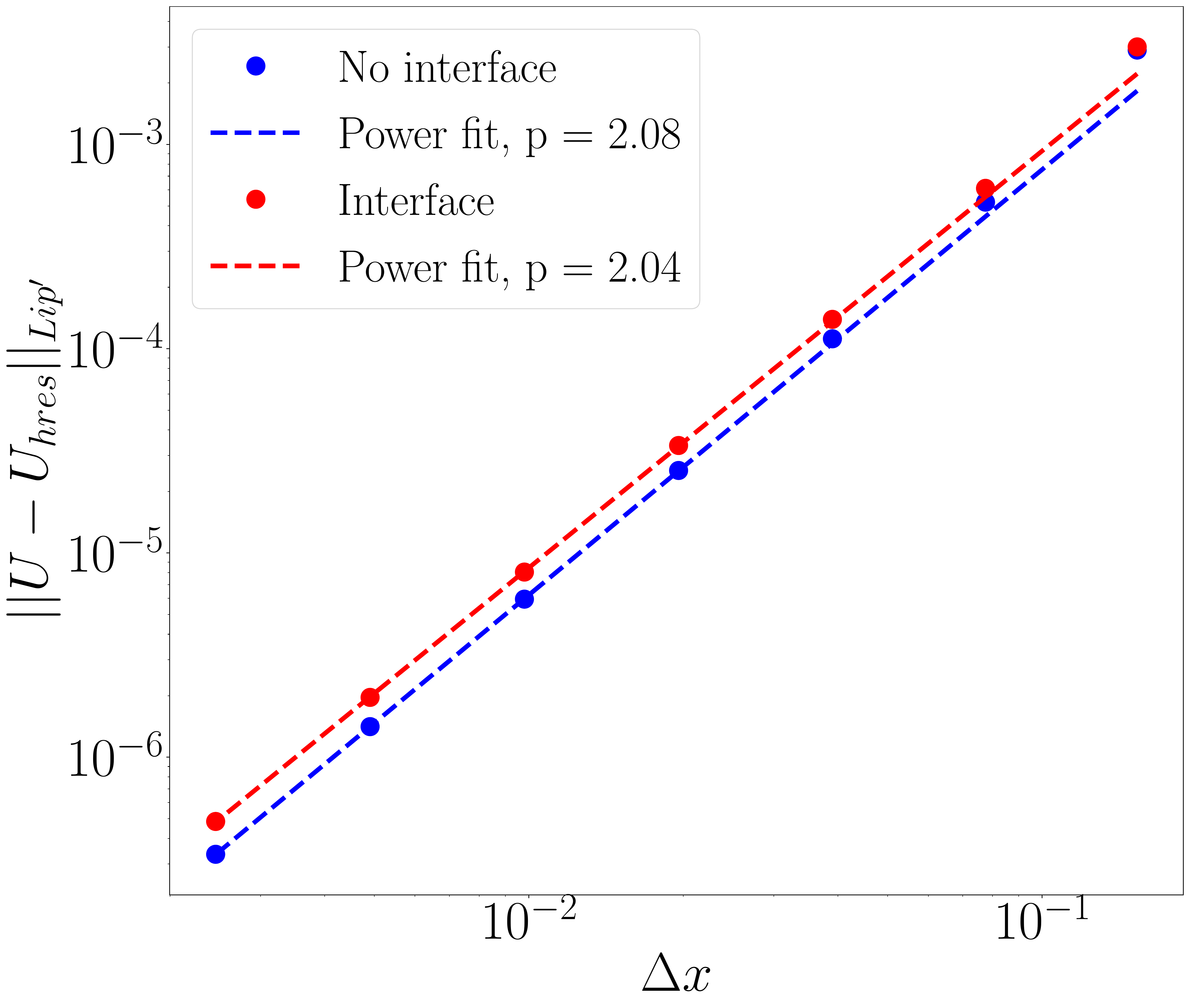}
    \caption{$Lip'$ convergence of \eqref{eq:sininidat} for the Burgers equation. $u_{hres}$ corresponds to the approximation with $N = 20480$.}
    \label{fig:sinconvergence2}
    \end{subfigure}
    \caption{Results for the one-dimensional Burgers problem at $T = 114$ and with $CFL = 0.1$.}
\end{figure}

We also continue the evolution up to $T = 114.0$, to allow the shock to travel through the interface several times. Given that the speed of the shock is 0.5, this means that it passes through the interface around $9$ times. Second-order convergence is maintained even after these 9 crosses, as can be seen in Figures \ref{fig:sinevolutionlong} and \ref{fig:sinconvergence2}.

\begin{figure}[H]
        \centering
        \includegraphics[width=0.8\textwidth]{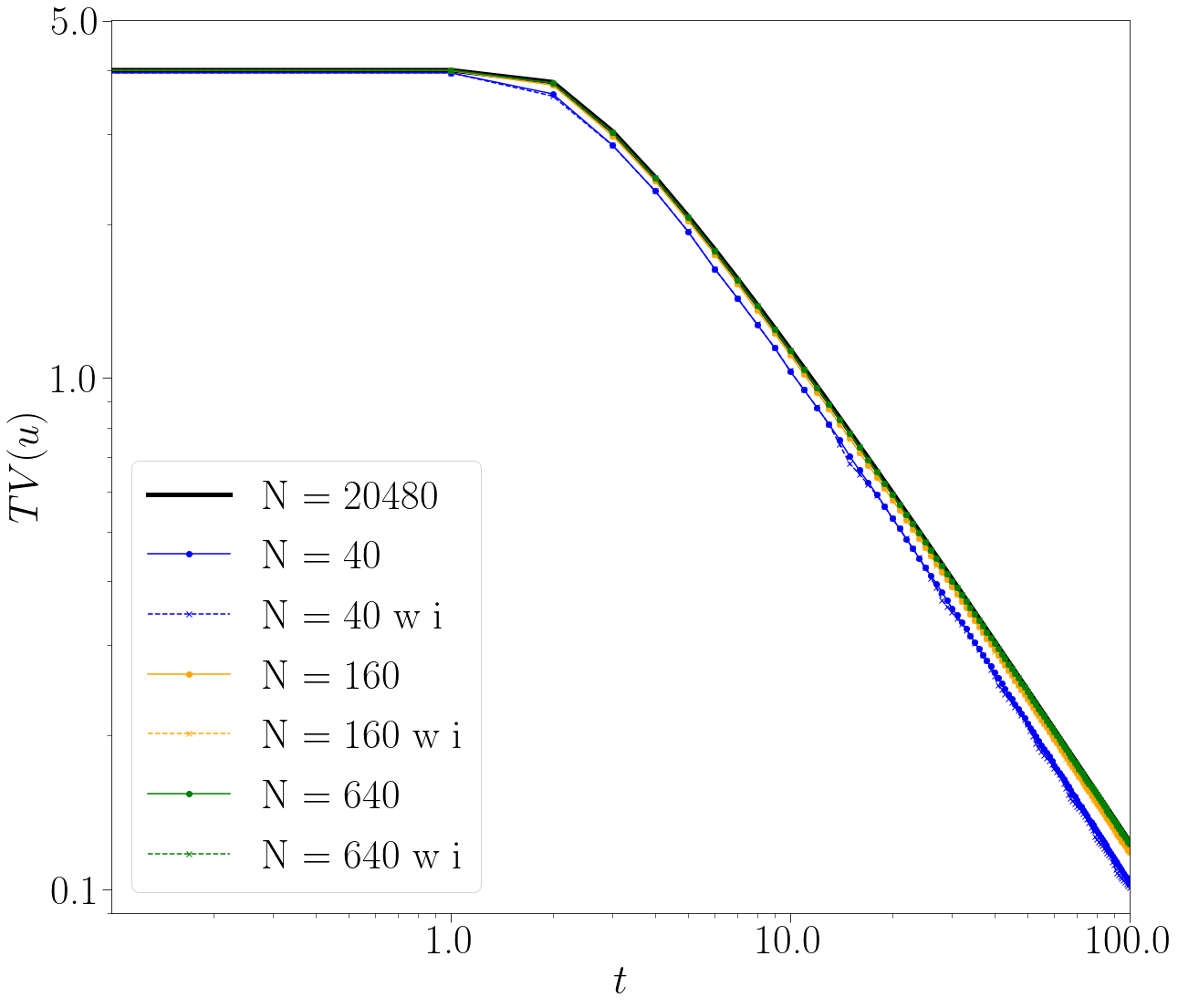}
    \caption{Logarithmic plot of the evolution of the total variation for the Burgers equation with (w i) and without an interface. The total variation is constant until the shock develops.}
    \label{fig:tv-burgers-sin}
\end{figure}

Figure \ref{fig:tv-burgers-sin} shows the total variation evolution of the solution of the Burgers equation until $T = 100$ for several grid resolutions. As expected, the total variation starts to diminish rapidly as soon as the shock develops. In contrast with the advection equation, the effect of the interface on the total variation is almost unnoticeable. Notice that for this case the TV plot is logarithmic.

\subsection{Two-dimensional Euler equation}

We evolve the two-dimensional Euler equations of gas dynamics \cite{LeVeque-1992},

\begin{equation}
    \frac{\partial}{\partial t}\begin{pmatrix}
\rho\\ 
\rho u\\
\rho w\\
E 
\end{pmatrix}
+ 
\frac{\partial}{\partial x}\begin{pmatrix}
\rho u\\ 
\rho u^{2} + p\\
\rho u w\\
u(E+p) 
\end{pmatrix}
+ 
\frac{\partial}{\partial y}\begin{pmatrix}
\rho u\\ 
\rho u w\\
\rho w^{2} + p\\
u(E+p) 
\end{pmatrix}
= 0,
\label{eq:euler2D}
\end{equation}

\noindent where $\rho$ is the density of the gas, $(u,w)$ the 2D fluid velocity, $E$ the energy density of the system and $p$ the pressure, which obeys state equation

\begin{equation}
    E = \frac{1}{2} \rho (u^{2} + w^{2}) + \frac{p}{\gamma -1},
    \label{eq:state-equation}
\end{equation}

\noindent with $\gamma$ the ratio of specific heats. We use $\gamma = 1.4$, which corresponds to an ideal gas.

\subsubsection{Implosion}

\begin{figure} [H]
    \centering
    \includegraphics[width = 0.4\linewidth]{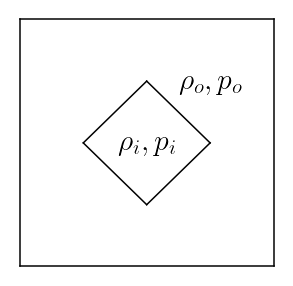}
    \caption{Initial data for the implosion problem.}
    \label{fig:implosion-inidat}
\end{figure}

We consider the implosion problem from \cite{liska2003}. It consists of a box $(x,y) \in [-0.3, 0.3]\times[-0.3,0.3]$ with reflective boundary conditions, with a smaller box inside rotated by $\pi/4$ with corners in points $[\pm 0.15, 0]$ and $[0, \pm0.15]$. The fluid velocity is set to zero everywhere, and initial pressure and density are set at $\rho_{i}=0.125$ and $p_{i} = 0.14$ in the inner box, and $\rho_{o} = 1, p_{o} = 0$ outside of it, as seen in figure \ref{fig:implosion-inidat}.

We take advantage of the symmetry of the problem and only evolve the upper right quadrant of the full domain, that is, $(x,y) \in (0,0.3)\times(0,0.3)$.

To impose the boundary conditions we add two extra rows of points outside of each boundary and set their field values to mirror the ones inside, except for the normal component of the momentum which is opposite to the mirror one inside the domain.


\squeezeup
\begin{figure}[H]
    \centering
    \includegraphics[width = 0.9\linewidth]{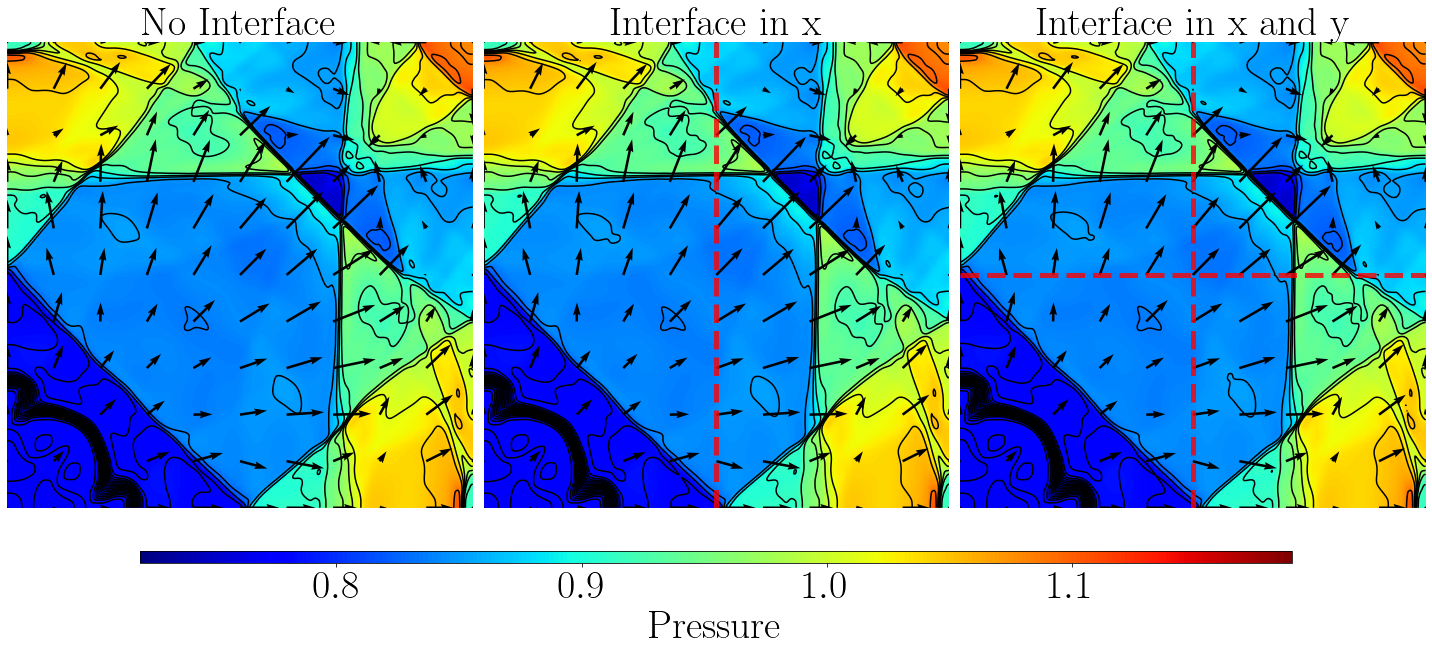}
    \caption{Implosion problem on a $400 \times 400$ grid with $CFL = 0.1$ at $T = 2.5$ with and without interfaces. Density contour lines and velocity arrows are overlaid over the pressure color maps.}
    \label{fig:implosion-t25}
\end{figure}
\squeezeup
\begin{figure}[H]
    \centering
    \includegraphics[width = 0.7\linewidth]{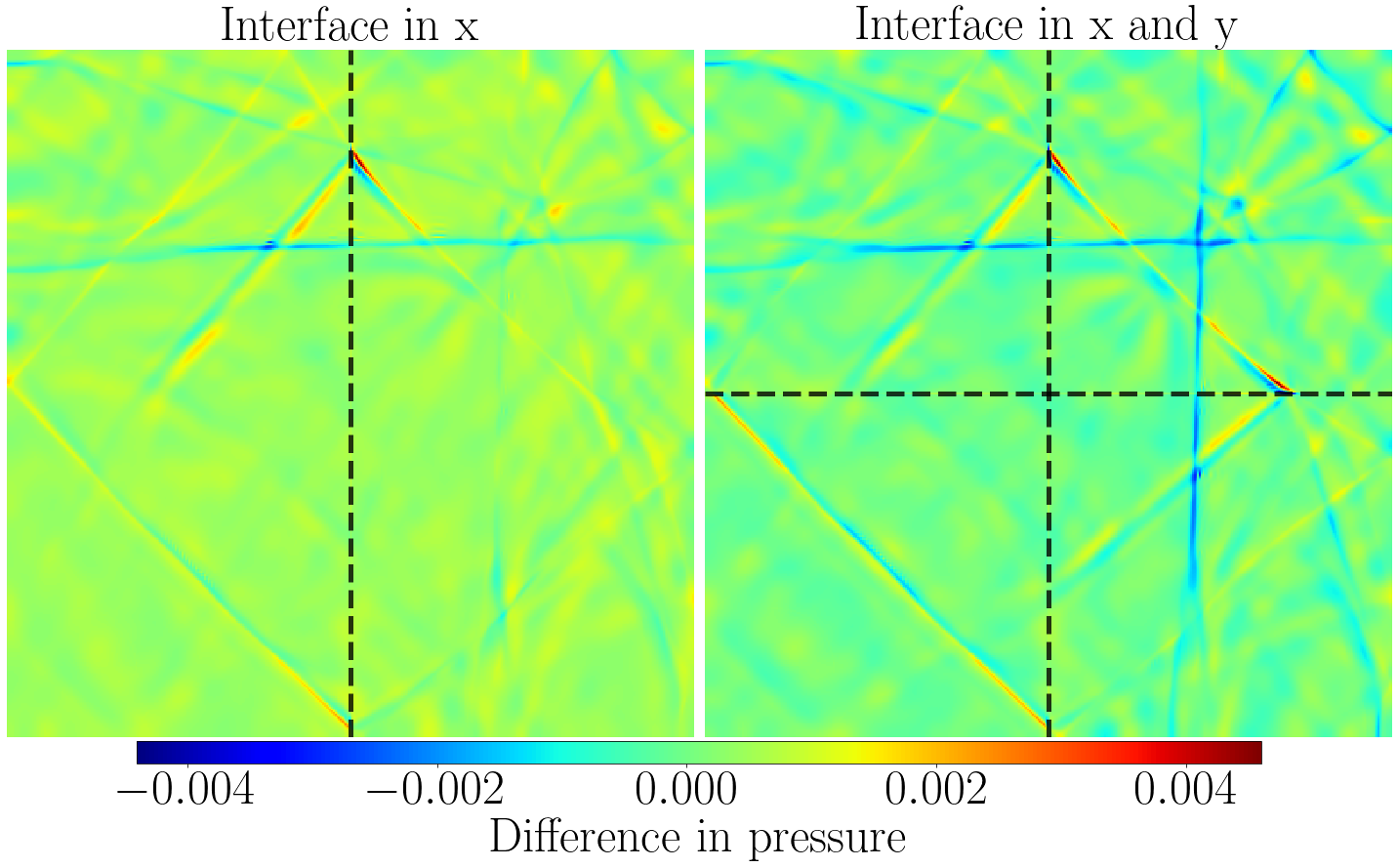}
    \caption{Difference between the approximations with and without interfaces for the implosion problem on a $400 \times 400$ grid with $CFL = 0.1$ at $T = 2.5$.}
    \label{fig:implosion-diff-t25}
\end{figure}
\squeezeup

We evolve this systems on $400\times400$ grids until $T = 2.5$ for three different cases: one with no interface, one with a vertical interface at $x = 0.15$, and the last one with two interfaces; one vertical at $x = 0.15$ and another horizontal at $y = 0.15$.
Results of the implementation of the method can be seen in figure \ref{fig:implosion-t25}. All three plots are very similar, which shows that the interface method works as intended. To help distinguish between them, we subtract the pressure approximations with interfaces from the pressure approximation without any interface and show them in Figure \ref{fig:implosion-diff-t25}.

We find that shocks are correctly propagated through the interfaces even when there is an angle between the propagation direction and the interface normal, and there is no indication of bouncing effects. Most of the differences between approximations are located around shocks. This is in full accordance with the differences seen in the one-dimensional case, although the $\theta$ value used here is different than that used in the scalar one-dimensional cases.

\subsubsection{Gresho vortex}

As a final test we evolve the Gresho vortex \cite{liska2003} in a box $(x,y) \in [-1,1)\times [-1,1)$ with periodic boundary conditions:

\begin{equation}
\renewcommand\arraystretch{1.3}
u_{\phi}(r) = \left\{ 
\begin{matrix*}[l]
5r \\ 2-5r  \\ 0  
\end{matrix*}\right. , \;
p(r) = \left\{ 
\begin{matrix*}[l]
5+\frac{25}{2}r^{2} & 0 \leq r < 0.2
\\ 9 +\frac{25}{2}r^{2}-20r+4\ln(5r)& 0.2 \leq r < 0.4 \\
3+4\ln(2) & 0.4 \leq r
\end{matrix*}\right.
\renewcommand\arraystretch{1}
\label{eq:gresho-inidat}
\end{equation}

\begin{figure}[H]
    \centering
        \centering
        \includegraphics[width = 0.8\linewidth]{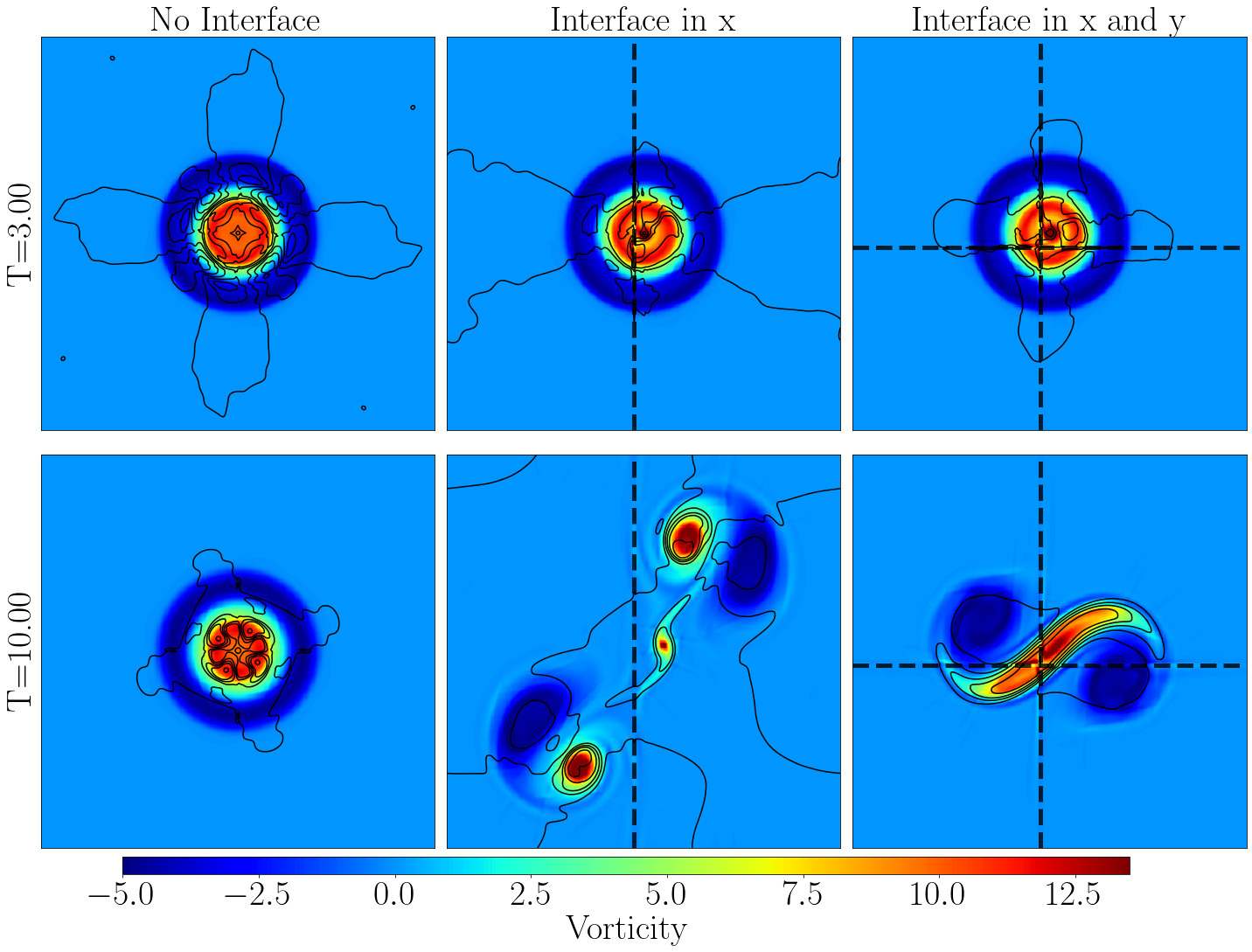}
    \caption{Gresho problem on a $400 \times 400$ grid with and without interfaces with $CFL = 0.1$. The upper row corresponds to $T = 3.0$ and the bottom row to $T = 10.0$. The interfaces are placed slightly off center, so as to break the symmetry of the problem. Density contours are overlaid over the vorticity color maps.}
    \label{fig:gresho-sol}
\end{figure}
The exact solution of this problem is independent of time, as the pressure gradient is constructed in such a way as to exactly cancel the centrifugal force given by the velocity profile. We solve this problem in a $400\times400$ grid with and without interfaces in the same way we did in the implosion problem but placing the interfaces slightly off-center so as to break the spatial symmetry of the numerical implementation. This symmetry breaking is realized by the interface points being treated differently than other points in the grid. This way we slightly excite unstable perturbations. Otherwise, if the interfaces are at the center of the vortex, the numerical simulation lasts a much longer time without unleashing the instabilities.

\begin{figure}[H]
    \centering
        \centering
        \includegraphics[width = 0.8\linewidth]{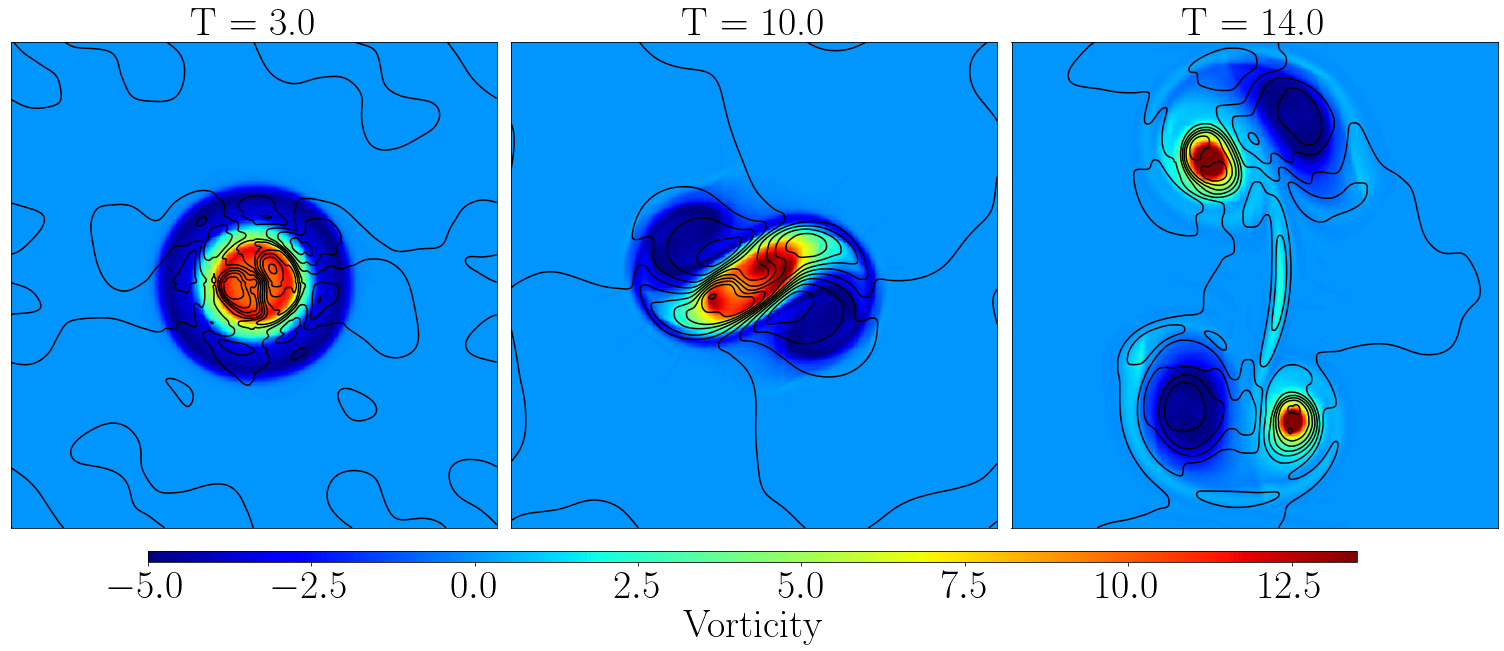}
    \caption{Gresho problem on a $400 \times 400$ grid with initial velocity perturbed by multiplying it by $(1.0 + 0.01\cos(\phi))$.}
    \label{fig:gresho-cosine-perturbation}
\end{figure}

Results seen in Figure \ref{fig:gresho-sol} show that the vortex is maintained for small times, but becomes unstable after a long evolution. This indicates that the usual bar instabilities of this problem are triggered much earlier due to the interfaces of this method. To check that we perturbed the initial data with a small azimuthal function and got similar effects. Figure \ref{fig:gresho-cosine-perturbation} shows the results of that perturbation, corroborating that the instability is both physical and generic. 

\newpage

\section{Conclusions:}

We introduced a new simple interface method for dealing with shock propagation in a multi-block setting. This method preserves the TVD condition while being second-order convergent. It has minimal communication among multiblocks, allowing for simple and efficient code parallelization. This method is based on the Kurganov-Tadmor semidiscrete scheme introduced in \cite{KT99}, and works by lowering the local order of convergence in the interface region to first order. As expected, the global convergence remains of second-order in the $Lip'$ norm, as was shown to be the case in all one-dimensional numerical experiments.
We also applied our method to the 2D Euler equations for two different situations: an implosion problem and a Gresho vortex. We found that the scheme allows the transmission of shock discontinuities without loss of phase or change in speed regardless of the orientation of the shock with respect to the interface normal. Finally, we also found that the bar instabilities of the Gresho vortex problem are triggered much earlier when interfaces not preserving the quadrant symmetry are present.

\section*{Acknowledgments}
The authors would like to thank funding from the Consejo Nacional de Investigaciones
Científicas y Técnicas (CONICET), Fondo para la Investigación Científica y Tecnológica (FONCYT), and Secretaría de Ciencia y Técnica, Universidad Nacional de Córdoba and discussions with Manuel Tiglio. 


\bibliography{mybib}
\bibliographystyle{ieeetr}

\end{document}